\documentclass{amsart}

\newcommand{\Z}{\ensuremath{\mathbb{Z}}}

\newcommand{\R}{\ensuremath{\mathbb{R}}}
\newcommand{\T}{\ensuremath{\mathbb{T}}}
\newcommand{\G}{\ensuremath{\mathcal{G}}}
\newcommand{\M}{\ensuremath{\mathcal{M}}}
\renewcommand{\L}{\ensuremath{\mathcal{L}}}

\newcommand{\spsi}{{\mbox{${\rm supp}~\hat\psi$}}} 
\newcommand{\mnp}{{\mbox{$2^{n-m-2}(\xi+2\cdot 2^m\pi)$}}}
\newcommand{\mnm}{{\mbox{$2^{n-m-2}(\xi-2\cdot 2^m\pi)$}}}

\newcommand{\ol}{\overline}
\newcommand{\be}{\begin{equation}}
\newcommand{\ee}{\end{equation}}
\newcommand{\bea}{\begin{eqnarray}}
\newcommand{\eea}{\end{eqnarray}}
\newcommand{\bes}{\begin{eqnarray*}}
\newcommand{\ees}{\end{eqnarray*}}

\newtheorem{theorem}{Theorem}[section]
\newtheorem{lemma}[theorem]{Lemma}
\newtheorem{proposition}[theorem]{Proposition}

\theoremstyle{definition}
\newtheorem{definition}[theorem]{Definition}
\newtheorem{remark}[theorem]{Remark}

\allowdisplaybreaks[4]
\numberwithin{equation}{section}

\begin{document}

\title{A Class of non-MRA Band-limited Wavelets}

\author{Biswaranjan Behera}
\address{Department of Mathematics, Indian Institute of Technology, 
Kanpur 208016, India}
\email{bbehera@iitk.ac.in}
\thanks{The first author was supported by The National Board 
for Higher Mathematics (NBHM), Govt. of India.}

\author{Shobha Madan}
\address{Department of Mathematics, Indian Institute of Technology, 
Kanpur 208016, India}
\email{madan@iitk.ac.in}

\date{July 2, 2001.}
\keywords{wavelet, MRA-wavelet, multiresolution analysis, 
dimension function}

\begin{abstract}
We give a characterization of a class of band-limited    
wavelets of $L^2({\mathbb R})$ and show that none of 
these wavelets come from a multiresolution analysis 
(MRA). For each $n\geq 2$, we construct a subset $S_n$ 
of ${\mathbb R}$ which is symmetric with respect to 
the origin. We give necessary and sufficient conditions 
on a function $\psi\in L^2({\mathbb R})$ with 
supp $\hat\psi\subseteq S_n$ to be an orthonormal 
wavelet. This result generalizes the characterization 
of a class of wavelets of E. Hern\'andez and G. Weiss. 
The dimension functions associated with these wavelets 
are also computed explicitly. Starting from the 
wavelets we have constructed, we are able to construct 
examples of wavelets in each of the equivalence 
classes of wavelets defined by E. Weber.
\end{abstract}

\maketitle

\section{Introduction} 
A function $\psi\in L^2(\R)$ is said to be an 
{\it orthonormal wavelet}, or simply a {\it wavelet}, if 
the system of functions $\{ \psi_{j,k}:j,k\in \Z \}$ 
forms an orthonormal basis for $L^2(\R)$, where
\[ 
\psi_{j,k}(x)=2^{j/2}\psi(2^jx-k),\quad  j,k\in \Z. 
\]
A function $\psi\in L^2(\R)$ is a wavelet of $L^2(\R)$ 
if and only if $\psi$ satisfies the following four conditions:
\begin{align}
& \sum_{j\in \Z} |\hat \psi (2^j\xi )|^2 
= 1\quad{\rm for~ a.e.}~\xi\in\R. \label{eqn:1W1} \\
& \sum_{j\geq 0} \hat \psi (2^j\xi)
\overline {\hat \psi \left(2^j(\xi +2m\pi)\right)}
= 0\quad{\rm for~ a.e.}~\xi\in\R~
{\rm and~ for~all}~m\in 2\Z + 1. \label{eqn:1W2} \\
& \sum_{k\in \Z}|\hat \psi (\xi +2k\pi)|^2 = 1
\quad{\rm for~a.e.}~\xi\in\R. \label{eqn:1W3} \\
& \sum_{k\in \Z} \hat\psi (\xi +2k\pi) 
\overline{\hat\psi\left(2^j (\xi + 2k\pi)\right)} = 0
\quad{\rm for~a.e.}~\xi\in\R~{\rm and~for~ all}~ j\geq 1.
\label{eqn:1W4}
\end{align}                 
We use the following definition of the Fourier transform
\[ 
\hat f(\xi)=\int_{\R}f(x)e^{-i\xi x}dx,\quad\xi\in\R. 
\]

The fact that equations (\ref{eqn:1W1})--(\ref{eqn:1W4}) 
characterize all wavelets of $L^2(\R)$ was observed by many 
authors. Lemari\'{e} and Meyer \cite{lm2} obtained these 
equations for compactly supported wavelets. Bonami, Soria 
and Weiss \cite{bsw} proved this fact under the assumption 
that $\psi$ is band-limited, i.e., $\hat\psi$ is compactly 
supported. For a function $\psi$ to be a wavelet, it is 
necessary that $\|\psi\|_2=1$. With this assumption on $\psi$, 
the following theorem was proved independently by 
G. Gripenberg \cite{grip} and X. Wang \cite{wan} 
(see also \cite{hkls}).

\begin{theorem}\label{thm:wavelet1}
Let $\psi\in L^2(\R)$ with $\|\psi\|_2=1$. Then $\psi$ is 
a wavelet of $L^2(\R)$ if and only if $\psi$ satisfies 
{\rm (\ref{eqn:1W1})} and {\rm(\ref{eqn:1W2})}.
\end{theorem}

For a proof of this theorem we refer to Chapter 7 of 
\cite{hw}. It was observed in \cite{bsw} that the 
orthonormality of the system $\{\psi_{j,k}:j,k\in\Z\}$ is 
characterized by (\ref{eqn:1W3}) and (\ref{eqn:1W4}), and 
completeness by (\ref{eqn:1W1}) and (\ref{eqn:1W2}) 
(see also \cite{hw}). A consequence of 
\mbox{Theorem \ref{thm:wavelet1}} is that equations 
(\ref{eqn:1W1}) and (\ref{eqn:1W2}), along with the 
assumption $\|\psi\|_2=1$, imply the equations 
(\ref{eqn:1W3}) and (\ref{eqn:1W4}).

If $\psi$ is a wavelet, then the support of $\hat\psi$ 
must have measure at least $2\pi$. This minimal measure 
is achieved only when $|\hat\psi|$ is the characteristic 
function of some measurable subset $K$ of $\R$. Such a 
wavelet is called a {\it minimally supported frequency} 
(MSF) wavelet and the associated set $K$ is called a 
{\it wavelet set}. We refer to \cite {hw} for proofs of 
the above statements.

A function is said to be {\it band-limited} if its Fourier 
transform has compact support. The simplest example of a 
band-limited wavelet is the Shannon wavelet whose Fourier 
transform is the characteristic function of the set 
$S=[-2\pi,-\pi]\cup [\pi,2\pi]$. Lemari\'{e} and Meyer 
\mbox{(\cite{lm}, \cite{mey1})} constructed a band-limited 
wavelet belonging to the Schwartz class. This wavelet 
satisfies $\hat\psi(\xi)=e^{i\xi/2}b(\xi)$, where $b$ is 
an even, non-negative ``bell-shaped" function and its support 
is equal to the set $[-\frac{8}{3}\pi,-\frac{2}{3}\pi] 
\cup [\frac{2}{3}\pi,\frac{8}{3}\pi]$. In \cite{bsw} 
the authors characterized all wavelets $\psi$ such that 
$|\hat\psi|$ is an even continuous function whose support 
is equal to the set 
$[-2\pi-\epsilon',-\pi+ \epsilon]\cup[\pi-\epsilon,2\pi+\epsilon']$, 
where $\epsilon,\epsilon'>0$ and $\epsilon + \epsilon' \leq \pi$, 
so that this class includes the Lemari\'{e}-Meyer wavelet as a 
particular case. Later, in \cite{hww} (see also \mbox{Theorem 4.1}, 
\mbox{Chapter 3} of \cite{hw}), the restriction on $\psi$ 
(i.e., $|\hat\psi|$ is even and continuous) was removed and all 
wavelets with Fourier transform supported in 
$[-\frac{8}{3}\pi,-\frac{2}{3}\pi] \cup [\frac{2}{3}\pi,\frac{8}{3}\pi]$ 
were characterized\footnote{In fact, they characterized all 
wavelets with Fourier transform supported in the set 
\mbox{$[-\frac{8}{3}a,4\pi-\frac{4}{3}a]$}, \mbox{$0<a\leq\pi$}. 
It was also shown that for such a wavelet $\psi$, it is 
necessary that $\hat\psi=0$ a.e. 
on $[-\frac{2}{3}a,2\pi-\frac{4}{3}a]$.}.

\medskip

This article is organised as follows. In \S 2 we construct 
a set $S_n$ for each $n\geq 2$, which is symmetric with 
respect to the origin and is a union of four intervals. 
In \S 3 all wavelets with Fourier transform supported in 
$S_n$ are characterized. Next, in \S 4 we show that none of 
these wavelets is associated with an MRA. For a wavelet to 
be associated with an MRA, it is necessary and sufficient that 
the corresponding dimension function is equal to 1 a.e. 
In \S 5 we compute the dimension function explicitly and show 
that it is an even function and takes all integral values 
from 0 to the maximum value it attains. This provides an 
alternative proof of the fact that these wavelets are not 
MRA-wavelets. An equivalence relation on the set of all 
wavelets of $L^2(\R)$ was defined in \cite{web}. Let $\M_n$ 
be the corresponding equivalence classes. Using the wavelets 
characterized in \S 3, we construct in \S 6 examples of wavelets 
in each $\M_n,n\geq 1$. Finally, we construct a family of wavelets 
in the equivalence class $\M_0$.

\section{Construction of the set $S_{n}$}
Let $n\geq 2$. Put 
\[ 
\begin{array}{ll}
  a_{n} = \frac{2^{n-1}}{2^n-1}\pi,        
& b_n = 2a_{n} = \frac{2^n}{2^n-1}\pi, \\
  c_{n} = \frac{2^{n-1}(2^n-2)}{2^n-1}\pi, 
& d_{n} = 2^n a_n = \frac{2^{2n-1}}{2^n-1}\pi. 
\end{array} 
\]
Define $S_{n} = S_{n}^{+} \cup S_{n}^{-}$, 
where $S_{n}^{+} = [a_{n},b_{n}]\cup [c_{n},d_{n}]$ 
and $S_n^- = -(S_n^+)$.

\begin{remark}
\begin{enumerate}
\item If $n = 2$, then $a_{2} = \frac{2}{3}\pi,
\ b_{2}=c_{2} = \frac{4}{3}\pi$
and $d_{2} = \frac{8}{3}\pi$. So $S_{2}$ is the set associated 
with the Lemari\'{e}-Meyer wavelets.
\item As $n$ increases, the set $[a_n,b_n]$ moves closer 
to $0$ whereas $[c_n,d_n]$ moves farther away. Also observe 
that the measure of the set $S_n$ is 
$2\bigl(\frac{2^{n-1}}{2^n-1}\pi +\frac{2^n}{2^n-1}\pi\bigr)$ 
which approaches $3\pi$ as $n$ tends to infinity.
\end{enumerate}
\end{remark}  

All terms in the equations (\ref{eqn:1W1})--(\ref{eqn:1W4}) 
that characterize all wavelets of $L^2(\R)$ involve dilations 
by powers of 2 and translations by integral multiples of $2\pi$. 
In the following lemma, for various subsets of $S_n$, we 
identify those translates and dilates which are possibly 
in $S_n$. We present the lemma in a tabular form for easy 
reference.

\begin{lemma}\label{lem:table}
For $n\geq 3$, let $S_{n}$ be as above, and let 
$e_{n} = \frac{2^n-2}{2^n-1}\pi$. Then
\vskip .2cm
\begin{center}
\begin{tabular}{|c|l|l|}    \hline   
$\xi\in the~set$&$\xi+2k\pi\not\in S_n$ 
unless & $2^j\xi\not\in S_n$ 
unless  \\ \hline \hline
$[a_{n},e_{n}]$ & $k = 0$               
                & $j = 0$             \\ \hline 
$[e_{n},b_{n}]$ & $k = 0$ or $-1$       
                & $j = 0$ or $n-1$    \\ \hline
$[c_{n},d_{n}]$ & $k = 0$ or $-2^{n-1}$ 
                & $j = 0$ or $-(n-1)$ \\ \hline
$[-e_n,-a_n]  $ & $k = 0$               
                & $j = 0$             \\ \hline
$[-b_n,-e_n]  $ & $k = 0$ or $1$        
                & $j = 0$ or $n-1$    \\ \hline
$[-d_n,-c_n]  $ & $k = 0$ or $2^{n-1}$  
                & $j = 0$ or $-(n-1)$ \\ \hline 
\end{tabular}
\end{center}
\end{lemma}
\vskip .2cm

Observe that 
\begin{align*}
[e_{n},b_{n}]-2\pi & =[-b_{n},-e_{n}], \\ 
 2^{n-1}[e_{n},b_{n}] & =[c_{n},d_{n}], \\ 
\intertext{and} 
[c_n,d_n]-2\cdot 2^{n-1}\pi & = [-d_{n},-c_{n}]. 
\end{align*}

\section{The characterization}
\label{sec:main}
The following theorem characterizes all wavelets $\psi$  
such that $\spsi$ is contained in $S_n$. The case $n=2$ is 
Theorem 4.1 of Chapter 3 in \cite{hw} in which case the 
wavelets are MRA-wavelets.

\begin{theorem}\label{thm:sn}
Let $~n\geq 2,~ \psi\in L^2(\R)$, {\rm supp} 
$\hat\psi\subseteq S_{n}$ and $~ b(\xi) =| \hat\psi(\xi)|$. 
Then $\psi$ is a wavelet for $L^2(\R )$ if and only if
\begin{enumerate}
\item [(i)]   $b(\xi) = 1$  
\quad for a.e. $\xi\in[a_{n},e_{n}] \cup [-e_n,-a_n]$,
\item [(ii)]  $b^2(\xi) + b^2(2^{n-1}\xi) = 1$
\quad for a.e. $\xi\in [e_n,b_n]$, 
\item [(iii)] $b^2(\xi) + b^2(\xi-2\pi) = 1  $
\quad for a.e. $\xi\in[e_n,b_n]$,
\item [(iv)]  $b(\xi) = b\left(2^{n-1}(\xi-2\pi)\right)$
\quad for a.e. $\xi\in[e_n,b_n]$,
\item [(v)]   $\hat \psi(\xi) = e^{i\theta(\xi)}b(\xi),$
where $\theta$ satisfies
\[ 
\theta(\xi)+\theta\left(2^{n-1}(\xi-2\pi)\right)
-\theta(\xi-2\pi)-\theta(2^{n-1}\xi)
= \left(2m(\xi)+1\right)\pi, 
\]
for some $ m(\xi)\in \Z,~$ for a.e. 
$\xi \in [e_{n},b_{n}]\cap ({\rm supp}~b) \cap 
(\frac{1}{2^{n-1}} {\rm supp}~b).$
\end{enumerate}
Moreover, for $n\geq 3$, none of these wavelets arise from an MRA.
\end{theorem}

\begin{remark} 
\begin{enumerate}
\item Observe that it follows from Theorem \ref{thm:sn} that 
$|\hat\psi|$ is completely determined by its values on 
$[e_{n},b_{n}].$ On this set, let $b=|\hat\psi|$ be an 
arbitrary measurable function taking values between 0
and 1. Using the relations following Lemma \ref{lem:table}, 
$|\hat\psi|$ can be extented to other sets of $S_n$ with the 
help of properties (i)--(iv). Properties (ii), (iii) and (iv) 
can also be written as 
\begin{eqnarray} 
b^2(\xi)+b^2\Bigl(\frac{1}{2^{n-1}}\xi\Bigr) & = & 1 
\quad{\rm for~a.e.~}\xi\in[c_{n},d_{n}], \label{eqn:cd3} \\ 
b^2(\xi)+b^2(\xi +2\pi ) & = & 1 
\quad{\rm for~a.e.~}\xi\in[-b_{n},-e_{n}], \label{E.2} \\
b\left(\frac{1}{2^{n-1}}\xi+2\pi\right) & = & b(\xi) 
\quad{~\rm for~a.e.~} [-d_{n},-c_{n}]. \label{E.3}
\end{eqnarray} 
\item If $({\rm supp}~b)\cap(\frac{1}{2^{n-1}}{\rm supp}~b)$ 
has an empty interior in $[e_{n},b_{n}]$, then $\theta $ can 
be chosen to be any measurable function. In particular,
we can take $\theta(\xi) = 0.$
\end{enumerate}
\end{remark} 

\noindent{\bf Proof of Theorem \ref{thm:sn}}
 
First let us assume that $\psi$ is a wavelet for $L^2(\R)$, 
$\spsi \subseteq S_{n}$ and $b = |\hat\psi|$. So $\psi$ 
satisfies equations (\ref{eqn:1W1})--(\ref{eqn:1W4}).
Consider equality (\ref{eqn:1W3}) :
\[
\sum_{k \in \Z} b^2(\xi+2k\pi) = 1 \quad {\rm for~a.e.}~\xi\in\R. 
\]
We use Lemma \ref{lem:table} to pick out the non-zero terms 
in this sum. If $|\xi|\in [a_{n},e_{n}]$, then only $k = 0$ 
will contribute to the sum. So we get $b(\xi)=1$ a.e., which 
is (i) of the theorem. For $\xi \in [e_{n},b_{n}]$, we get 
non-zero contributions from $k = 0$ and $k = -1$. That is, 
$b^{2}(\xi ) + b^{2}(\xi -2\pi ) = 1$ for a.e. $\xi$ which 
proves (iii). Also if $\xi\in [e_{n},b_{n}]$, then 
$2^{n-1}\xi \in [c_{n},d_{n}]$, and we get (non-zero terms 
now correspond to $k = 0$ and $k = -2^{n-1}$)    
\begin{equation}\label{eqn:cd1}
b^2(2^{n-1}\xi) + b^2(2^{n-1}\xi-2^n\pi) = 1 
\quad{\rm for~a.e.~}\xi\in[e_n,b_n]. 
\end{equation}
This can also be written as
\begin{equation}\label{eqn:cd2}
b^2(\xi ) + b^2(\xi - 2^n\pi ) = 1
\quad{\rm for~a.e.~}\xi\in[c_{n},d_{n}]. 
\end{equation}
Now consider equality (\ref{eqn:1W1}):
\[
\sum_{j\in \Z} b^{2}(2^{j}\xi ) = 1
\quad{\rm for~a.e.}~\xi\in\R. 
\]
For $\xi\in[e_{n},b_{n}]$, the only possible non-zero terms 
correspond to $j = 0$ and $j=n-1$. So we conclude 
$b^{2}(\xi )+b^{2}(2^{n-1} \xi ) = 1$. This proves (ii). 
Combining (ii) and (\ref{eqn:cd1}) we get condition (iv). 
It remains to prove condition (v). For this purpose 
we consider equality (\ref{eqn:1W4}) with $j = n-1$:
\[
\sum_{k\in\Z} \hat\psi(\xi+2k\pi)
\overline{\hat\psi\left(2^{n-1}(\xi + 2k\pi)\right)} = 0  
\quad {\rm for~ a.e.}~ \xi\in\R.
\]
If $\xi \in [e_{n},b_{n}]$, the non-zero terms in the sum 
come from $k = 0$ and $-1$. Hence,
\begin{equation}\label{eqn:eb}
\hat \psi (\xi )\overline{\hat \psi (2^{n-1}\xi )} +
\hat \psi (\xi -2\pi ) 
\overline{\hat \psi \left(2^{n-1}(\xi -2\pi )\right)} =0 
   \quad{\rm for~a.e.}~\xi\in [e_{n},b_{n}].
\end{equation}
Using (iii), (\ref{eqn:cd1}) and (\ref{eqn:eb}), we see that 
for almost every $\xi\in [e_{n},b_{n}]$, the vectors 
\[
\left(\hat \psi (\xi ),\hat \psi (\xi -2\pi )\right)
\quad{\rm and}\quad
\left(\hat\psi(2^{n-1}\xi),\hat\psi(2^{n-1}(\xi-2\pi))\right) 
\]
are orthonormal in ${\mathbb C}^2$. If we let $\hat\psi(\xi) 
= e^{i\theta(\xi)}b(\xi)$, then it follows that, for some 
real-valued measurable function $\alpha$, 
\begin{eqnarray}
\lefteqn{e^{i\alpha (\xi )}\left( e^{i\theta (\xi )} b(\xi ) ,
     e^{i\theta (\xi-2\pi )}b(\xi -2\pi )\right)} \nonumber \\
 & & = \left(-e^{-i\theta (2^{n-1}
         (\xi -2\pi))}b(2^{n-1}(\xi -2\pi)) ,
         e^{-i\theta(2^{n-1}\xi )}b(2^{n-1}\xi )\right)
\label{eqn:phase} 
\end{eqnarray}
for a.e. $\xi\in [e_n,b_n]$. But from (ii), (iii) and (iv), 
we know that
\begin{equation}
 \left. 
\begin{array}{rcl}
 b(2^{n-1}\xi ) & = & b(\xi - 2\pi)  \\
 b(\xi )        & = & b\left( 2^{n-1}(\xi - 2\pi )\right)
\end{array} \right\}\quad\mbox{for ~a.e}~\xi\in [e_n,b_n]. 
\label{eqn:bxi} 
\end{equation}
So (\ref{eqn:phase}) can be written as
\[
 \left. \begin{array}{rcl}
  \left[e^{i\alpha (\xi )}e^{i\theta (\xi )} +
 e^{-i\theta \left(2^{n-1}(\xi -2\pi )\right)}\right]b(\xi ) & = & 0  \cr
  \left[e^{i\alpha (\xi )}e^{i\theta (\xi - 2\pi )} - 
 e^{-i\theta (2^{n-1}\xi )}\right]b(2^{n-1}\xi ) & = & 0
\end{array} \right\} \quad\mbox{for~a.e.}~\xi\in [e_n,b_n]. 
\] 
This shows that
\[ 
\begin{array}{rcll}
e^{-i\alpha (\xi )} & = & e^{i[\theta (\xi ) + 
\theta (2^{n-1}(\xi - 2\pi)) + \pi]}
 & {\rm for~a.e.}~\xi\in[e_n,b_n] \cap {\rm supp}~b. \\
e^{-i\alpha (\xi )} & = & e^{i[\theta (\xi - 2\pi ) + 
\theta (2^{n-1}\xi )]}
 & {\rm for~a.e.}~\xi\in[e_n,b_n] \cap 
(\frac{1}{2^{n-1}}{\rm supp}~b). 
\end{array} 
\]
Hence, for almost every $\xi\in [e_n,b_n]
\cap {\rm supp}~b \cap(\frac{1}{2^{n-1}} {\rm supp}~b)$, 
we have
\[
\theta(\xi)+\theta\left(2^{n-1}(\xi -2\pi)\right)-
\theta(\xi-2\pi)-\theta (2^{n-1}\xi ) 
= \left(2m(\xi )+1\right)\pi, 
\]
for some integer-valued measurable function $m$. This proves (v).

{\bf We now prove the converse}. Suppose $\psi\in L^{2}(\R)$, 
$\spsi\subseteq S_{n}$ and the function $b(\xi) = |\hat\psi(\xi )|$ 
satisfies conditions (i)--(v) of the theorem.  
By Theorem \ref{thm:wavelet1}, to show that $\psi$ is a 
wavelet, it is sufficient to 
show that $\|\psi\|_2 = 1$ and  $\psi$ satisfies 
(\ref {eqn:1W1}) and (\ref{eqn:1W2}). We have,
\begin{eqnarray*}
  2\pi\|\psi\|^{2}_{2}  &  = & \|\hat\psi\|^{2}_{2} =
   \int_{S_n} |\hat\psi (\xi )|^2~ d\xi \\
  & = & \left( \int_{-d_n}^{-c_n} + \int_{-b_n}^{-e_n} 
+ \int_{-e_n}^{-a_n}+ \int_{a_n}^{e_n} + \int_{e_n}^{b_n} 
+ \int_{c_n}^{d_n}\right) b^2(\xi )~d\xi \\
  & = & \int_{-d_n}^{-c_n} + \int_{c_n}^{d_n} + 
\int_{-b_n}^{-e_n} + \int_{e_n}^{b_n} + 2(e_n-a_n).
\end{eqnarray*}
By changing variables $\xi\rightarrow 2^{n-1}(\xi-2\pi)$, 
$\xi\rightarrow 2^{n-1}\xi$ and $\xi\rightarrow\xi-2\pi$ 
in the first, second and third integrals respectively, we get
\begin{eqnarray*} 
 2\pi\|\psi\|^2_2 & = & 2^{n-1}\int_{e_n}^{b_n}\left[b^2
              \left(2^{n-1}(\xi -2\pi )\right) 
               + b^2(2^{n-1}\xi )\right] d\xi \\
  &   &      + \int_{e_n}^{b_n}\left[b^2(\xi -2\pi)
             + b^2(\xi )\right] d\xi + 2(e_n-a_n) \\
 & = & 2^{n-1}(b_n-e_n) + (b_n-e_n) + 2(e_n-a_n) = 2\pi,
\end{eqnarray*}
where we have used (\ref{eqn:cd1}) and property (iii) of 
the theorem. Hence, $\|\psi\|_2 = 1$.

We will now show that $\psi$ satisfies (\ref{eqn:1W1}). Let
\[ 
\rho(\xi )=\sum_{j\in \Z}\mid\hat\psi(2^{j}\xi)\mid^{2}. 
\] 
Suppose $\xi >0 $. Observe that
\begin{eqnarray*} 
 \R^{+} & = & \bigcup_{l\in \Z} 2^{l} [a_{n},b_{n}] \\ 
        & = & \bigcup_{l \in \Z} 2^{l} \left([a_{n},e_{n}]
\cup [e_{n},b_{n}]\right). 
\end{eqnarray*}
So there is an $l\in\Z$ such that 
$\xi\in 2^l[a_n,e_n]\cup 2^l[e_n,b_n]$. 
If $\xi\in 2^l[a_n,e_n]$ then since $2^{-l}\xi\in [a_n,e_n]$, 
by Lemma \ref{lem:table} we have, $2^j(2^{-l}\xi )\not\in S_n$ 
if $j\neq 0$. That is, $2^j\xi\not\in S_n$ if $j\neq -l$. 
Hence, $\rho (\xi ) = 1$, by (i). Similarly by using (ii), we 
can prove that $\rho (\xi ) = 1$ if $\xi\in 2^l[e_n,b_n]$. 
A similar decomposition for $\xi < 0$ proves that 
$\rho (\xi )$ = 1 for a.e. $\xi\in\R$.

Finally, we have to show that $\psi$ satisfies (\ref{eqn:1W2}). 
For $m\in 2\Z+1$, let us denote the function on the left hand 
side of (\ref{eqn:1W2}) by $t_m(\xi)$. Then
\begin{eqnarray}
t_{m}(\xi ) & = & \sum_{j\geq 0} \hat \psi (2^{j}\xi)
 \overline{\hat\psi \left(2^j(\xi +2m\pi)\right)} \nonumber \\
            & = & \sum_{j\geq 0} \hat \psi 
\left(2^{j}(\xi+2m\pi-2m\pi)\right)
 \overline{\hat\psi \left(2^j(\xi +2m\pi)\right)} \nonumber \\
            & = & \ol{t_{-m}(\xi+2m\pi)} \label{eqn:tm}.
\end{eqnarray}
Therefore, we have only to show that $t_{m}(\xi ) = 0$ for 
a.e. $\xi$, if  $m\in 2\Z_{-}+1$ (negative odd integers).
Suppose $m\in 2\Z_{-}+1$ and $m\neq -1$. Let $\xi \in \R$ 
and suppose that $2^{j}\xi \in S_{n}$. Since $j\geq$ 0, 
$n\geq 2$ and $m\neq -1$ and is odd, therefore  
$2^j m\neq 0,\pm 1,
\pm 2^{n-1}$. So by Lemma \ref{lem:table}, we have 
$2^j\xi + 2\cdot 2^jm\pi \not\in S_n$. This implies that 
each term of $t_{m}(\xi )$ is zero which proves that 
$t_{m}$ is zero. It now remains to show that $t_{-1}(\xi ) = 0$ 
for a.e $\xi$. We have
\[ 
t_{-1}(\xi ) = \sum_{j\geq 0} \hat\psi (2^{j}\xi )
\overline{\hat\psi (2^{j}\xi -2\cdot 2^{j}\pi )}. 
\]
Let $\xi\in\R$ and suppose that $2^j\xi\in S_{n}$. Then 
again by Lemma \ref{lem:table}, 
$2^{j}\xi -2\cdot 2^{j}\pi \not\in S_{n}$ if 
$-2^{j}\neq 0, \pm 1, \pm 2^{n-1}.$
So the only possible $j$'s to contribute a non-zero 
term are $j = 0$ and $j = n-1$. Thus,
\begin{equation}
 t_{-1}(\xi ) = \hat\psi (\xi )\overline{\hat\psi (\xi -2\pi )}
 +\hat\psi(2^{n-1}\xi )
\overline{\hat\psi (2^{n-1}\xi -2\cdot 2^{n-1}\pi )}.  
\label{eqn:tq}
\end{equation}
Now, both $\xi$ and $\xi -2\pi$ belong to $S_{n}$ only if 
$\xi \in[e_{n},b_{n}]$. Hence, the first term of (\ref{eqn:tq}) 
is zero unless $\xi \in   [e_{n},b_{n}]$. Similarly, both 
$2^{n-1}\xi$ and $2^{n-1}\xi -2\cdot2^{n-1}\pi$ belong to $S_{n}$ 
only when $2^{n-1}\xi \in [c_{n},d_{n}],$ which is equivalent 
to saying that $\xi \in [e_{n},b_{n}]$. That is, the second 
term of (\ref{eqn:tq}) is also zero unless
$\xi \in [e_{n},b_{n}]$. Thus we get, $t_{-1}(\xi ) = 0$ 
if $\xi\not\in  [e_{n},b_{n}]$. Now on $[e_{n},b_{n}]$, 
we have, by (\ref{eqn:bxi})
\[ 
t_{-1}(\xi)=b(\xi)b(2^{n-1}\xi) e^{i[\theta (\xi)- 
\theta (\xi -2\pi)]} + b(2^{n-1}\xi)b(\xi)
e^{i[\theta (2^{n-1}\xi )+ 
\theta \left(2^{n-1}(\xi -2\pi)\right)]}. 
\]
If $\xi\not\in [e_{n},b_{n}] \cap {\rm supp}~b 
\cap (\frac{1}{2^{n-1}}{\rm supp}~b)$,
then either $b(\xi ) = 0$ or $b(2^{n-1}\xi ) = 0$. 
So $t_{-1}(\xi ) = 0$.
And if $\xi \in [e_{n},b_{n}] \cap {\rm supp}~b 
\cap(\frac{1}{2^{n-1}}{\rm supp}~b)$, then by (v), 
$t_{-1}(\xi ) = 0$. 
This completes the characterization.
\qed
\vskip .2cm

There is another equation which characterizes all wavelets 
$\psi$ such that $\hat\psi$ is supported in the set $S_n$ 
and the function $b=|\hat\psi|$ is even. For $n=2$, the 
following proposition is proved in \cite{hw} (see 
Proposition 4.7, Chapter 3). We observe that 
the result can be extended to the general case.

\begin{proposition}
Suppose that $\psi$ is a wavelet of $L^2(\R)$, $b=|\hat\psi |$, 
and ${\rm supp}~b\subseteq S_n$. Then $b$ is almost everywhere 
even if and only if
\begin{equation} \label{E5}
b^2(\xi)+b^2(2\pi-\xi)=1\quad for~a.e.~ \xi\in [e_n,b_n].
\end{equation}
\end{proposition}

\noindent{\bf Proof:}
Let $\psi$ be a wavelet of $L^2(\R)$ such that $b=|\hat\psi|$ 
is supported in $S_n$. Suppose that $b$ is an even function 
and $\xi\in [e_n,b_n]$. Since $-\xi\in [-b_n,-e_n]$, 
we have
\begin{eqnarray*} 
1 &=& b^2(-\xi)+b^2(-\xi+2\pi),\quad {\rm by~(\ref{E.2})}\\
&=& b^2(\xi)+b^2(2\pi-\xi),\quad {\rm since}~ b ~{\rm is~even},
\end{eqnarray*}
which is (\ref{E5}).

Conversely, suppose that (\ref{E5}) holds. Since by (i) of 
Theorem \ref{thm:sn}, $b$ is even on $[a_n,e_n]$, it is 
enough to show that $b$ is even on the sets $[e_n,b_n]$ 
and $[c_n,d_n]$.

\noindent(a) Let $\xi\in [e_n,b_n]$. Therefore, 
$-\xi\in [-b_n,-e_n]$. Then
\begin{equation*} 
b^2(-\xi)+b^2(-\xi+2\pi)=1,\quad {\rm by~(\ref{E.2}).}
\end{equation*}
This fact, together with (\ref{E5}), gives us
\begin{equation} \label{E6}
b(\xi)= b(-\xi).
\end{equation}
(b) Now let  $\xi\in [c_n,d_n]$ . Therefore, 
$\frac{1}{2^{n-1}}\xi\in [e_n,b_n]$. 
 From (\ref{E6}) we get
\begin{equation} \label{E7}
b\Bigl(\frac{1}{2^{n-1}}\xi\Bigr) 
= b\Bigl(- \frac{1}{2^{n-1}}\xi \Bigr).
\end{equation}
Since $- \frac{1}{2^{n-1}}\xi \in[-b_n,-e_n]$, 
using (\ref{E.2}) we obtain
\begin{equation} \label{E8}
b^2\Bigl(-\frac{1}{2^{n-1}}\xi\Bigr )
+b^2\Bigl(- \frac{1}{2^{n-1}}\xi +2\pi\Bigr)=1.
\end{equation}
 Now, since $-\xi\in[-d_n,-c_n]$, we get (using (\ref{E.3}))
\begin{equation} \label{E9}
b\Bigl(- \frac{1}{2^{n-1}}\xi +2\pi\Bigr)=b(-\xi).
\end{equation}
Substituting (\ref{E7}) and (\ref{E9}) in (\ref{E8}), we get
\begin{equation} \label{E10}
b^2\Bigl(\frac{1}{2^{n-1}}\xi \Bigr)+b^2(-\xi)=1.
\end{equation}
Comparing (\ref{E10}) and (\ref{eqn:cd3}) we get 
$b(\xi)=b(-\xi)$. This proves the proposition.
\qed

\section{The wavelets associated with $S_n$, $n\geq 3$ are non-MRA}

A multiresolution analysis (MRA) is a sequence  of closed 
subspaces \mbox{$\{V_j:j\in \Z\}$} of $L^2(\R)$ satisfying 
the following properties:
\begin{enumerate}
\item[(i)] $V_j\subset V_{j+1}$ for all $j\in\Z$
\item[(ii)] $\cup_{j\in\Z}V_j$ is dense in $L^2(\R)$ and 
$\cap_{j\in\Z}V_j=\{0\}$
\item[(iii)] $f\in V_j$ if and only if $f(2\cdot)\in V_{j+1}$ 
for all $j\in\Z$
\item[(iv)] there exists a function $\varphi \in L^2(\R)$ 
such that $\{\varphi(\cdot-k):k\in\Z\}$ forms an orthonormal 
basis for $V_0$.
\end{enumerate}
A function $\varphi$ that satisfies property (iv) is called a 
{\it scaling function} for the MRA. If $\varphi$ is a scaling 
function for a given MRA, then it is easy to see that 
(see \cite{hw}) there exists a $2\pi$-periodic function $m_0$ 
in $L^2(\T)$, called the {\it low-pass filter}, such that
\begin{equation}
\hat\varphi(2\xi)=m_0(\xi)\hat\varphi(\xi).
\label{eqn:lpf}
\end{equation}
The function $m_0$ satisfies the equation 
\begin{equation}\label{eqn:m0xi}
|m_0(\xi)|^2+|m_0(\xi+\pi)|^2=1\quad{\rm for~ a.e}~\xi\in \T.
\end{equation} 
Using (\ref{eqn:lpf}) and (\ref{eqn:m0xi}), one can then 
construct a wavelet $\psi$ associated with the MRA. For any 
such wavelet the following equation holds (see \cite{hw} 
for details):
\begin{equation}
|\hat\varphi(\xi)|^2=\sum\limits_{j\geq 1}|\hat\psi(2^j\xi)|^2. 
\label{eqn:sf}
\end{equation} 
We will now prove the last statement in Theorem \ref{thm:sn}.

\begin{proposition}
If $n\geq 3$, then the wavelets characterized in 
Theorem {\rm \ref{thm:sn}} are not associated with any MRA.
\end{proposition}

\noindent{\bf Proof:} 
Let $\psi$ be any wavelet such that $\hat\psi$ is supported 
in $S_n$ and let it be associated with an MRA. Let $\varphi$ 
be the corresponding scaling function and $m_0$ be the low-pass 
filter associated with $\varphi$. Using (\ref{eqn:sf}) we can 
easily find that 
\begin{eqnarray}
|\hat\varphi(\xi)|=
\begin{cases}
1               & {\rm if}~~|\xi|\leq a_n\\
b(2^{n-l-1}\xi) & {\rm if}~~|\xi|\in 2^l[e_n,b_n],
                  ~~0\leq l\leq n-2\\
0               &{\rm otherwise.}
\end{cases}
\label{eqn:phihat1}
\end{eqnarray}
Hence,
\begin{eqnarray}
|\hat\varphi(2\xi)|=
\begin{cases}
1             & {\rm if}~~|\xi|\leq \frac{a_n}{2}\\
b(2^{n-l}\xi) & {\rm if}~~|\xi|\in 2^{l-1}[e_n,b_n],
                 ~~0\leq l\leq n-2\\
0             &{\rm otherwise.}
\end{cases}
\label{eqn:phihat2}
\end{eqnarray}
{\bf Case 1:} $b(\xi)\not\equiv 1$ on $[e_n,b_n]$,  i.e., 
$b(\xi)\not\equiv 0$ on $[c_n,d_n]$.
\vskip .2cm

If $|\xi|\leq a_n$, then from (\ref{eqn:lpf}) and 
(\ref{eqn:phihat1}), we have 
$|\hat\varphi(2\xi)|=|m_0(\xi)|\cdot |\hat\varphi(\xi)|=|m_0(\xi)|$. 
Therefore, by (\ref{eqn:phihat2}), we obtain  
$|m_0(\xi)|= 1~{\rm if}~ |\xi|\leq \frac{a_n}{2}$.
Since $m_0$ is $2\pi$-periodic, we have
\[ 
|m_0(\xi)|=1\quad{\rm on}~\bigl[-\tfrac{a_n}{2},
\tfrac{a_n}{2}\bigr]+2\cdot 2^{n-3}\pi. 
\]
Note that $[-\frac{a_n}{2},\frac{a_n}{2}]+2
\cdot 2^{n-3}\pi =[\frac{c_n}{2},\frac{d_n}{2}]
=2^{n-2}[e_n,b_n]$. Now, from (\ref{eqn:phihat1}), 
on $ 2^{n-2}[e_n,b_n]$ we have $|\hat\varphi(\xi)|=b(2\xi)$. 
So $|\hat\varphi(2\xi)|=|m_0(\xi)|\cdot |\hat\varphi(\xi)|
=b(2\xi)$. 
But by (\ref{eqn:phihat2}), $|\hat\varphi(2\xi)|=0$ on 
$2^{n-2}[e_n,b_n]$. Therefore, $b(2\xi)=0$ on 
$2^{n-2}[e_n,b_n]$. That is, $b(\xi)=0$ on 
$2^{n-1} [e_n,b_n] = [c_n,d_n]$, which is a contradiction.
\vskip .2cm
\noindent{\bf Case 2:} $b(\xi)\equiv 1$ on $[e_n,b_n]$.
\vskip .2cm
As above, we have
\begin{eqnarray}
|\hat\varphi(\xi)|=
\begin{cases}
1   & {\rm if}~~|\xi|\leq a_n\\
    & {\rm or~if}~~\xi\in 2^l[-b_n,-e_n],~~0\leq l\leq n-2\\
0   & {\rm otherwise.}
\end{cases}
\label{eqn:phihat3}
\end{eqnarray}
Therefore,
\begin{eqnarray}
|\hat\varphi(2\xi)|=
\begin{cases}
1   & {\rm if}~~|\xi|\leq \frac{a_n}{2}\\
    & {\rm or~if}~~\xi\in 2^{l-1}[-b_n,-e_n],
      ~~0\leq l\leq n-2\\
0   & {\rm otherwise.}
\end{cases}
\label{eqn:phihat4}
\end{eqnarray}
Now $|\hat\varphi(2\xi)|=|m_0(\xi)|\cdot|\hat\varphi(\xi)|
= |m_0(\xi)|$ if $|\xi|\leq a_n$, 
by (\ref{eqn:phihat3}). Therefore, by (\ref{eqn:phihat4}), 
$|m_0(\xi)|=1$ if $|\xi|\leq \frac{a_n}{2}$. Using the 
$2\pi$-periodicity of $m_0$, we get
\[
|m_0(\xi)|=1 \quad{\rm on}~~\bigl[-\tfrac{a_n}{2},
\tfrac{a_n}{2}\bigr]-2\cdot 2^{n-3}\pi. 
\]
Note that  
$[-\frac{a_n}{2},\frac{a_n}{2}]-2\cdot 2^{n-3}\pi
=[-\frac{d_n}{2},-\frac{c_n}{2}] 
= 2^{n-2}[-b_n,-e_n]$. Hence, on the interval 
$2^{n-2}[-b_n,-e_n]$, we have 
\[
|\hat\varphi(2\xi)|=|m_0(\xi)|\cdot|\hat\varphi(\xi)|= 
|m_0(\xi)|=1. 
\]
But by (\ref{eqn:phihat4}), $|\hat\varphi(2\xi)|=0$ on 
$2^{n-2}[-b_n,-e_n]$, which is again a contradiction. 
Therefore, the wavelets characterized in Theorem \ref{thm:sn} 
are non-MRA wavelets, if $n\ge 3$. 
\qed

\section{The dimension functions of wavelets associated with $S_n$}

Let $\psi$ be any wavelet of $L^2(\R)$. Define $V_j$ to be 
the closure of the span of $\{\psi_{l,k}:l<j,\ k\in\Z\}$. 
If $\{V_j:j\in\Z\}$ forms an MRA of $L^2(\R)$, then we say 
that $\psi$ is associated with an MRA, or $\psi$ is an 
MRA-wavelet.
 
Not every wavelet of $L^2(\R)$ is associated with an MRA. 
In the last section we proved that the wavelets associated 
with the set $S_n$ are non-MRA wavelets. The first example 
of a non-MRA wavelet, which is an MSF wavelet, was given by 
\mbox{J.L. Journ\'e}. The associated wavelet set is 
$J=\left[-\tfrac{32}{7}\pi,-4\pi\right]\cup 
\left[-\pi,-\tfrac{4}{7}\pi\right]
\cup \left[\tfrac{4}{7}\pi,\pi\right]\cup 
\left[4\pi,\tfrac{32}{7}\pi\right]$. 
Another interesting non-MRA wavelet is the Lemari\'{e} 
wavelet which is also an MSF wavelet 
with the corresponding wavelet set
$L=\left[-\tfrac{8}{7}\pi,-\tfrac{4}{7}\pi\right]\cup
\left[\tfrac{4}{7}\pi,\tfrac{6}{7}\pi\right] \cup
\left[\tfrac{24}{7}\pi,\tfrac{32}{7}\pi\right]$. 
These two wavelets belong to the class of wavelets 
characterized in \S \ref{sec:main}
(see Remark \ref{rem:other} below).

Given a wavelet $\psi$ of $L^2(\R)$, there is an associated 
function $D_\psi$, called the {\it dimension function}, 
defined by
\[ D_{\psi}(\xi ) = \sum_{j\geq 1}\sum_{k\in\Z}|\hat\psi 
                   \left(2^j(\xi + 2k\pi)\right)|^2. \]
It was proved independently by Gripenberg and Wang that a 
wavelet $\psi$ is an MRA-wavelet if and only if 
$D_{\psi}$ = 1 a.e. For a proof of this fact see \cite{hw}.

In this section we will compute the dimension functions for 
the wavelets characterized in Theorem \ref{thm:sn}. This will 
provide an alternative proof of the fact that these wavelets 
are not associated with any MRA. These dimension functions 
are symmetric with respect to the origin and attain all 
integral values starting from  $0$ to the maximum value. 
Also note that if $n\geq 3$, then all wavelets whose Fourier 
transform is supported in $S_n$ have the same dimension 
function $D_n$. We will prove the following proposition.

\begin{proposition}
Fix $n\geq 3$. Let $\psi$ be a wavelet such that 
$\spsi\subseteq S_n$ and let $D_n$ be the 
dimension function associated with $\psi$. Then 
\begin{eqnarray*}  
D_n(\xi) = 
  \begin{cases}
      n-1 & \text{a.e. if $|\xi|\in 
            \bigl[0,\frac{2}{2^n-1}\pi\bigr]$} \\ 
      r-1 & \text{a.e. if $|\xi|
            \in \bigl[\frac{2^{n-r}}{2^n-1}\pi,
            \frac{2^{n-r+1}}{2^n-1}\pi\bigr]
             ~~(2\leq r\leq n-1)$} \\  
      0 &  \text{a.e. if $|\xi|\in 
           \bigl[\frac{2^{n-1}}{2^n-1}\pi,
           \frac{2^n-2}{2^n-1}\pi\bigr] = [a_{n},e_{n}]$} \\    
      1 &  \text{a.e. if $|\xi|\in 
           \bigl[\frac{2^n-2}{2^n-1}\pi,\pi\bigr]
               = [e_{n},\pi ]$}.
   \end{cases}
\end{eqnarray*}
\end{proposition}

\noindent{\bf Proof:} For $l\in\Z$, we define 
\[ 
p_l = \frac{2^{n+l-1}}{2^n-1}\pi = 2^la_n\quad {\rm and}\quad 
q_l = \frac{2^{n+l-1}-2^l}{2^n-1}\pi = 2^{l-1}e_n. 
\]

Note that $p_0 = a_n$, $p_1 = 2p_0 = b_n$, $p_n = d_n$, 
$q_1 = e_n$ and $q_n = c_n.$ Also observe that 
$p_l<q_{l+1}<p_{l+1}$ for all $l\geq 0$. As $D_n$ is 
$2\pi$-periodic, it is enough to compute its values for 
$\xi\in [-\pi,\pi]$. Since $[-\pi,\pi]\subset[-p_1,p_1]$,
we will compute $D_n$ for the interval $[-p_1,p_1]$. 
Note that
\[ 
(0,p_1] = \bigcup_{l\geq 0}2^{-l}[p_0,p_1] 
        = \bigcup_{l\geq 0}[p_{-l},p_{-l+1}].  
\]
{\bf Case 1: $\xi\in[p_0,p_1]$}
\vskip .2cm

If $2^m\leq k\leq 2^{m+1}-1$ $(0\leq m\leq n-3)$, 
then by an elementary calculation we see that 
$\xi+2k\pi\in[p_{m+2},q_{m+3}]$. So $2^j(\xi+2k\pi)
\in[p_{j+m+2},q_{j+m+3}]$. If $j+m+3\leq n$ 
then $2^j(\xi+2k\pi)\in[p_1,q_n]=[b_n,c_n]$ which 
is not in $\spsi$, and if $j+m+3\geq n+1$ then 
$2^j(\xi+2k\pi)\geq d_n$, which shows that 
$2^j(\xi+2k\pi)\not\in\spsi$. Now if $k\geq 2^{n-2}$, 
then $\xi+2k\pi\geq d_n$. Thus, we have proved that if 
$k\geq 1$, then $2^j(\xi+2k\pi)\not\in\spsi$ for all 
$j\geq 1$. Similarly, we can show that if $k\leq -2$, 
then $2^j(\xi+2k\pi)\not\in\spsi$ for all 
$j\geq 1$. Therefore, we have only to consider $k=0,-1$. 
Note that $[p_0,p_1]=[p_0,q_1]\cup[q_1,p_1]$.
\vskip .2cm
(a) If $\xi\in[p_0,q_1],$ then $2^j\xi\in[p_j,q_{j+1}]$. 
Now $j\geq n\Rightarrow 2^j\xi\geq p_n=d_n$ and 
$j\leq n-1\Rightarrow j+1\leq n\Rightarrow$ for all 
$j\geq 1,~ 2^j\xi\in[p_j,q_n]\subset[p_1,q_n]=[b_n,c_n]$ 
which is not in $\spsi$. 
Also, $\xi-2\pi\in[-q_2,-p_1]\Rightarrow 2^j(\xi-2\pi)
\in[-q_{j+2},-p_{j+1}]$. 
A similar argument as in the case of $2^j\xi$ shows that 
$2^j(\xi-2\pi)\not\in\spsi$ for all $j\geq 1$. So   
$D_n(\xi)=0$ for a.e. $\xi\in[p_0,q_1]$.
\vskip .2cm
(b) If $\xi\in[q_1,p_1],$ then $2^j\xi\in[q_{j+1},p_{j+1}]$.  
So $2^j\xi\not\in\spsi$ if $j\not =n-1$. Similarly, 
$2^j(\xi-2\pi)\in[-p_{j+1},-q_{j+1}]$. 
So $2^j(\xi-2\pi)\not\in\spsi$ if $j\not =n-1$. Therefore, 
$D_n(\xi)=|\hat\psi (2^{n-1}\xi)|^2 + |\hat\psi 
\left(2^{n-1}(\xi - 2\pi)\right)|^2 =1$, 
by (\ref{eqn:cd2}), as $2^{n-1}\xi\in[q_n,p_n]=[c_n,d_n]$.
\vskip .2cm
We now proceed to compute $D_n$ on other sets. Observe that
\begin{eqnarray*}
D_n(\xi ) & = & \sum_{k=-\infty}^{-1}\sum_{j\geq 1}
\bigl|\hat\psi \left(2^j(\xi + 2k\pi)\right)
\bigr|^2 + \sum_{j\geq 1}\bigl|\hat\psi (2^j\xi)\bigr|^2  \\ 
&   & +\sum_{k=1}^{\infty}\sum_{j\geq 1}\bigl|\hat\psi 
\left(2^j(\xi + 2k\pi )\right)\bigr|^2 \\
& = & D_n^-(\xi )+D_n^0(\xi ) +D_n^+(\xi ), {\rm ~say}. 
\end{eqnarray*} 
Let $\xi\in\bigl[\frac{2^{n-l}}{2^n-1}\pi,
\frac{2^{n-l+1}}{2^n-1}\pi\bigr]
  =[p_{-l+1},p_{-l+2}]$, $l\geq 2$.
\vskip .2cm
As in the above cases, we can show that if 
$k\leq -2^{n-2}$, or if $-2^{m+1}-1\leq k\leq -2^m-1$  
$(0\leq m\leq n-3)$, then $2^j(\xi+2k\pi)\not\in\spsi$ 
for all $j\geq 1$. Also, if $k=-2^m (0\leq m\leq n-3)$, 
then $2^j(\xi+2k\pi)\not\in\spsi$ for all $j\geq 1$ 
unless $j=n-m-2$.
Hence, in the sum $D_n^-(\xi )$, we have only to consider 
$j=n-m-2$ $(0\leq m\leq n-3)$ 
and $k=-2^m$. Thus,
\[ 
D_n^-(\xi )=\sum_{m=0}^{n-3}\bigl|\hat\psi 
\left(2^{n-m-2}(\xi - 2\cdot 2^m\pi)\right)\bigr|^2.
\]
In a similar manner, we can show that
\[ 
D_n^+(\xi )=\sum_{m=0}^{n-3}\bigl|\hat\psi 
\left(2^{n-m-2}(\xi + 2\cdot 2^m\pi)\right)\bigr|^2. 
\]

\noindent{\bf Case 2:}
$\xi\in\bigl[\frac{2^{n-l}}{2^n-1}\pi,
\frac{2^{n-l+1}}{2^n-1}\pi\bigr]$, $l\geq n$ 
\vskip .2cm
By a straightforward calculation, one can show that if 
$l\geq n$, then 
\[
\mnp\in[c_n,d_n],\quad{\rm for}~~0\leq m\leq n-3.
\]
Therefore, 
\bes
&   & \hat\psi\bigl(\mnp\bigr)+\hat\psi\bigl(\mnm\bigr) \\
& = & \hat\psi\bigl(\mnp\bigr)+
      \hat\psi\bigl(\mnp-2^{n-1}\pi\bigr) \\
& = & 1,\quad{\rm by}~ (\ref{eqn:cd2}).
\ees
Hence, we get $D_n^+(\xi)+D_n^-(\xi)=n-2$.

Now, $[p_{-l+1},p_{-l+2}]=[p_{-l+1},q_{-l+2}]\cup 
[q_{-l+2},p_{-l+2}]$. 

(a) If $\xi\in [p_{-l+1},q_{-l+2}]$, then $2^j\xi\not\in\spsi$, 
if $j\not = l-1$. But $2^{l-1}\xi\in[p_0,q_1]=[a_1,e_n]$. 
Therefore, $D_n^0(\xi )=|\hat\psi(2^{l-1}\xi)|^2 =1$.

(b) If $\xi\in [q_{-l+2},p_{-l+2}]$, then $2^j\xi\not\in\spsi$ 
if $j\not = l-1,n+l-2$. Now $2^{n+l-2}\xi\in[q_n,p_n]=
[c_n,d_n]$. By (\ref{eqn:cd3}) we get
$D_n^0(\xi)=|\hat\psi(2^{n+l-2}\xi)|^2+ |\hat\psi(2^{l-1}\xi)|^2=1$.

Thus, in either case $D_n(\xi)=n-2+1=n-1$. We have proved that 
$D_n(\xi)=n-1$, 
if $\xi\in \Bigl[\frac{2^{n-l}}{2^n-1}\pi,
\frac{2^{n-l+1}}{2^n-1}\pi\Bigr],~l\geq n$. 
That is, 
$D_n(\xi)=n-1$ for a.e. $\xi\in[0,\frac{2}{2^n-1}\pi]$.
\vskip .2cm

\noindent{\bf Case 3:} $\xi\in\bigl[\frac{2^{n-l}}{2^n-1}\pi,
\frac{2^{n-l+1}}{2^n-1}\pi\bigr]$, 
$l=2$ 
\vskip .2cm
As in Case 2, it can be shown that $\mnm$ and $\mnp$ are not 
in $\spsi$ for all $m$ such that $0\leq m\leq n-3$. Therefore, 
$D_n^-(\xi)=D_n^+(\xi)=0$.

Now if $\xi\in\bigl[\frac{2^{n-2}}{2^n-1}\pi,
\frac{2^{n-1}-1}{2^n-1}\pi\bigr]$, then 
$2^j\xi\not\in\spsi$ if $j\not=1$; and if $\xi$ belongs 
to the set  
$\bigl[\frac{2^{n-1}-1}{2^n-1}\pi,\frac{2^{n-1}}{2^n-1}\pi\bigr]$, 
then $2^j\xi\not\in\spsi$ if $j\not=1,n$. In either case,
\[ 
D_n(\xi)=D_n^0(\xi)=1.
\]

\noindent{\bf Case 4:} $\xi\in\bigl[\frac{2^{n-l}}{2^n-1}\pi,
\frac{2^{n-l+1}}{2^n-1}\pi\bigr]$, 
$3\leq l\leq n-1$ (This case is required if $n\geq 4$.) 
\vskip .2cm
We can write $2^{n-m-2}\xi\in\bigl[\frac{2^{n+p-l}}{2^n-1}\pi,
\frac{2^{n+p}}{2^n-1}\pi\bigr]$, where $p=n-m-l-1.$ Using the 
condition on support of $\hat\psi,$ we can show that if 
$\xi\in\bigl[\frac{2^{n+p-l}}{2^n-1}\pi,
\frac{2^{n+p}}{2^n-1}\pi\bigr]$, then 
$\xi\pm 2\cdot 2^{n-2}\pi\in\spsi$ only when $p\leq -1$. 
So in order that 
$\hat\psi(2^{n-m-2}\xi\pm 2\cdot 2^{n-2}\pi)\not= 0$, 
we must have $p=n-m-l-1\leq -1$, i.e., $m\geq n-l$. We also 
have, $0\leq m\leq n-3$. Therefore, $1\leq n-m-2\leq l-2$. Thus
\begin{eqnarray*}
D_n(\xi) & = & \sum\limits_{j\geq 1}\bigl|\hat\psi(2^j\xi)\bigr|^2\\
&   & +\sum_{m=0}^{n-3}\left\{\bigl|\hat\psi 
\left(2^{n-m-2}(\xi + 2\cdot 2^m\pi)\right)\bigr|^2             
+ \bigl|\hat\psi \left(2^{n-m-2}(\xi - 2\cdot 2^m\pi)
\right)\bigr|^2\right\}\\ 
& = & \sum\limits_{j\geq 1}\bigl|\hat\psi(2^j\xi)\bigr|^2+ 
\sum_{s=1}^{l-2}\left\{\bigl|\hat\psi 
(2^s\xi + 2\cdot 2^{n-2}\pi)\bigr|^2 + 
\bigl|\hat\psi (2^s\xi - 2\cdot 2^{n-2}\pi)\bigr|^2\right\}.
\end{eqnarray*}
Now, $2^s\xi+2\cdot 2^{n-2}\pi\in[c_n,d_n],~1\leq s\leq l-2$. 
By (\ref{eqn:cd2}), the second sum on the right hand side of 
the last equality above is equal to $l-2$. Also as in the 
previous cases, the first sum can be shown to be equal to $1$. 
Hence, $D_n(\xi)=l-1$.	
\qed  
\medskip
\begin{remark}\label{rem:other}
In \cite{hkls} MSF wavelets were considered where they were 
called unimodular wavelets. Let $\psi$ be an MSF wavelet and 
$K = \spsi$. Define $K^+ = K\cap [0,\infty)$ and 
$K^- = K\cap (-\infty,0]$. One of the results proved in 
\cite{hkls} is about wavelets $\psi$ such that 
$K^- = -K^+$ and $K^+$ consists of two disjoint intervals. 
They proved the following:

(A) Let $\hat\psi(\xi) = \mu(\xi)\chi_K(\xi),$ where 
$K = K^+\cup K^-, K^- = -K^+$
and $K^+ = [a_n,\pi]\cup[2^{n-1}\pi,d_n]$ with $|\mu(\xi)|=1.$ 
Then $\psi$ is a wavelet for $L^2(\R)$. Moreover, each MSF 
wavelet for which $K^-=-K^+$, and $K^+$ is a union of two 
disjoint intervals is of this form.

If we define $b = \chi_{[e_n,\pi]}$ on $[e_n,b_n]$, and extend 
it to $S_n$ by using (i)--(iv) of Theorem \ref{thm:sn}, then 
we get the wavelets described in (A). In particular, for $n=3$, 
the corresponding wavelet is the Journ\'{e} wavelet. 
The paper \cite{hkls} also contains the following result.

(B) Let $\hat\psi(\xi)=\mu(\xi)\chi_K(\xi)$ where
\begin{align*}
K^- & = \left[-2\left(1-\frac{2p+1}{2^n-1}\right)\pi,
             -\left(1-\frac{2p+1}{2^n-1}\right)\pi\right] \\
\intertext{and} 
K^+ & = \left[\frac{2(p+1)}{2^n-1}\pi,\frac{2(2p+1)}{2^n-1}\pi\right]
         \cup\left[\frac{2^n(2p+1)}{2^n-1}\pi,
          \frac{2^{n+1} (p+1)}{2^n-1}\pi\right]   
\end{align*}
for $n\geq 3,~ 1\leq p\leq 2^{n-1}-2$ and $|\mu (\xi )| = 1.$
Then $\psi$ is a wavelet for $L^2(\R)$. Moreover, each MSF 
wavelet for which $K^-$ is an interval and $K^+$ is the union 
of two disjoint intervals is of this form.

If we define $b(\xi) = 0$ on the interval $[e_n,b_n]$ and 
extend it to $S_n$ by using (i)--(iv) of 
\mbox{Theorem \ref{thm:sn}}, then we get
\[ b = \chi_{[-b_n,-a_n]\cup [a_n,e_n]\cup [c_n,d_n]} \]
which corresponds to the case $p = 2^{n-2}-1$ in (B). 
In particular, for $n = 3$, we get the Lemari\'{e} wavelet.
\end{remark}

\section{Equivalence classes of wavelets}
First of all, we recall the equivalence relation on the set 
of all wavelets of $L^2(\R)$ defined in \cite{web}. Let $\psi$ 
be a wavelet of $L^{2}(\R)$. Define $V_j = \overline{sp}\{\psi_{l,k}:
\mbox{$l<j$},\mbox{$k\in\Z$}\}$. Then, it is easy to verify that 
the subspaces $V_j$, $j\in\Z$ satisfy the following properties:
\begin{enumerate}
\item $V_j\subset V_{j+1}$ for all $j\in\Z$
\item $f\in V_j$ if and only if $f(2\cdot)\in V_{j+1}$ 
for all $j\in\Z$
\item $\cup_{j\in\Z} V_j$ is dense in $L^{2}(\R)$ and 
$\cap_{j\in\Z} = \{0\}$
\item $V_0$ is invariant under the group of translation by integers.
\end{enumerate} 
A natural question is the following: Are there other groups of 
translations under which $V_0$ 
is invariant?
For $\alpha\in\R$, let $T_\alpha$ denote the unitary operator 
$T_\alpha f = f(\cdot -\alpha)$. Consider the group of 
translations $\G_n = \{T_{\frac{m}{2^n}}:m\in\Z\}$ and the 
group $\G_{\infty}=\{T_{\alpha }:\alpha\in\R\}$. Let $\L_n$ be 
the collection of all wavelets such that the corresponding 
space $V_0$ is invariant under the group $\G_n$. Then we have
\[ 
\L_{0}\supset\L_{1}\supset\L_{2}\supset\cdots\supset\L_{n}
\supset\L_{n+1} \supset\cdots\supset\L_{\infty}. 
\]   
An equivalence relation can now be defined on the set of 
all wavelets. The equivalence classes are given by 
$\M_n = \L_n\setminus \L_{n+1}$, with $\M_\infty = \L_\infty$. 
Therefore, $\M_n$, $n\geq 0$ is the class of wavelets such 
that $V_0$ is invariant under the group $\G_n$ but not under 
$\G_{n+1}$.

In \cite{web} Weber proved the following facts. 
\begin{enumerate}
\item $\M_\infty=\cap_{n\geq 0}\L_n.$
\item $\M_\infty$ is precisely the collection of all MSF wavelets. 
	
For the classes $\M_n,n\geq 1$ the following characterization 
is given.
\item The equivalence class $\M_n$, $n\geq 1$ consists of all 
wavelets $\psi$ such that $\spsi$ is not partially 
self-similar with respect to  any odd multiple of 
$2^k \pi,\ k=1,2,\dots,n$, but is partially self-similar with 
respect to some odd multiple of $2^{n+1} \pi$. 
\end{enumerate}

\begin{definition} {\rm \cite{web}}
A set $E$ is said to be partially self-similar with respect 
to $\alpha\in \R$ if there is a set $F$ of positive measure 
such that both $F$ and $F+\alpha$ are  subsets of $E$.
\label{def:pss}
\end{definition} 

In the same paper, Weber produced examples of wavelets 
belonging to the first few equivalence classes $\M_n$, 
namely for $n=0,1,2$ and $3$. This motivated us to 
construct wavelets belonging to each equivalence class 
$\M_n$. After we constructed these wavelets we became 
aware of the article \cite{sw}, in which S. Schaffer and 
\mbox{E. Weber} constructed wavelets for the classes 
$\M_n,n\geq 1$ by using the method of ``operator interpolation'' 
(see \cite{dl}). In fact, they did this for other dilation 
factors as well. In this section we will construct wavelets 
belonging to each of the classes $\M_n,n\geq 0$ by a 
different method. Our approach is simpler than that of 
\cite{sw} in the sense that for each $n\geq 3$, we 
construct a function $\psi_n$ such that $\psi_n$ has the 
required properties to be in $\M_{n-2}$ as characterized 
by Weber. Then we show that $\psi_n$ is a wavelet. Further, 
we construct a family of wavelets belonging to the 
equivalence class $\M_0$.

Note that, in view of the characterization of $\M_n$, 
to prove that $\psi_n\in\M_{n-2}$, it is sufficient to 
show that $\psi_n$ satisfies the following two conditions.
\begin{enumerate}
\item[a.] (supp $\hat\psi_n)\cap({\rm supp}~\hat\psi_n+2^k q\pi) 
= \emptyset$ for all $q\in 2\Z+1$ and $k=1,2,\dots n-2$.
\item[b.] There exists a subset $E$ of supp $\hat\psi_n$ 
such that $E+2^{n-1} q\pi\subset {\rm supp}~\hat\psi_n$ 
for some $q\in 2\Z+1.$	
\end{enumerate}

\subsection{Construction of wavelets in $\M_n$, $n\geq 1$}
\label{subsec:equiv1}

Our starting point is the MSF wavelet $\gamma_n$, $n\geq 3$ 
given via its Fourier transform:
\[ 
\hat\gamma_n=\chi_{_{W_n}}, \quad{\rm where}\quad 
W_n=[-b_n,-a_n]\cup[a_n,e_n]\cup[c_n,d_n]. 
\]
The wavelet $\gamma_n$ belongs to the class of band-limited 
wavelets characterized in \mbox{Theorem \ref{thm:sn}} 
(see Remark \ref{rem:other}(B)).

We will construct $\psi_n$ from this function in the 
following manner. We translate the interval
$[\frac{a_n}{2},\frac{e_n}{2}]+2^{n-1}\pi$ (which is 
a subset of $[c_n,d_n]$) to the left by a factor of 
$2^{n-1}\pi$ and assign values $\frac{1}{\sqrt 2}$ to 
$\hat\psi_n$ on both these sets. Then, we translate 
$[a_n,e_n]$ to the right by a factor of $2^n\pi$ and 
assign to $\hat\psi_n$ the value $\frac{1}{\sqrt 2}$ on 
$[a_n,e_n]$ and $-\frac{1}{\sqrt 2}$ on $[a_n,e_n]+2^n\pi$. 
Assign the value $1$ to $\hat\psi_n$ on the remaining sets  
of $W_n$ and $0$ elsewhere. More precisely, we have the 
following function. 
\begin{eqnarray*}
\hat\psi_n(\xi)=
\begin{cases}
1  & {\rm if}~\xi\in[-b_n,-a_n]\cup[c_n,\frac{a_n}{2}+2^{n-1}\pi]
     \cup[\frac{e_n}{2}+2^{n-1}\pi,d_n] \\
\frac{1}{\sqrt 2}  & {\rm if}~\xi\in[\frac{a_n}{2},\frac{e_n}{2}] 
     \cup[a_n,e_n]\cup[\frac{a_n}{2}+2^{n-1}\pi,
                 \frac{e_n}{2}+2^{n-1}\pi] \\
-\frac{1}{\sqrt 2} & {\rm if}~\xi\in[a_n+2^n\pi,e_n+2^n\pi] \\
0                  & {\rm otherwise.}
\end{cases}
\end{eqnarray*}

Let us make some observations about the set 
$F_n={\rm supp}~\hat\psi_n.$ It is easy to 
verify the following facts, stated in a tabular form 
for easy reference.
\vskip .2cm
\begin{center}
\begin{tabular}{|c|l|l|}    \hline   
$\xi\in {\rm the ~set}$ & $\xi +2k\pi \not\in F_n$ 
unless  & $2^j\xi\not\in F_n$ unless  
  \\ \hline \hline
$[-b_n,-e_n]$                 &  $k = 0$            
& $j = 0$\\  \hline 
$[-e_n,-a_n]$                 &  $k = 0 $           
& $j = 0$     \\  \hline
$[\frac{a_n}{2},\frac{e_n}{2}]$ & $k = 0,2^{n-2}$ 
& $j = 0,1$ \\ \hline
$[a_{n},e_{n}]$                 &  $k = 0,2^{n-1}$    
& $j = 0,-1$   \\  \hline
$[c_{n},\frac{a_n}{2}+2^{n-1}\pi]$ & $k = 0$  
& $j = 0$ \\ \hline 
$[\frac{a_n}{2}+2^{n-1}\pi,\frac{e_n}{2}+2^{n-1}\pi]$    
                              &  $k = 0,-2^{n-2}$  
& $j = 0,1$ \\ \hline 
$[\frac{e_n}{2}+2^{n-1}\pi,d_{n}]$  &  $k = 0$ 
& $j = 0$ \\ \hline 
$[a_n+2^n\pi,e_n+2^n\pi]$ & $k = 0,-2^{n-1}$ 
& $j = 0, -1$ \\ \hline 
\end{tabular}
\end{center}
\vskip .2cm
\begin{theorem} 
For each $n\geq 3$, the function $\psi_n$ defined above 
is a wavelet and belongs to the equivalence class $\M_{n-2}$.
\end{theorem}
\noindent{\bf Proof:} 
To prove that $\psi_n$ is a wavelet, it is sufficient 
to show that $\psi$ satisfies the following three 
properties (see Theorem \ref{thm:wavelet1}).
\begin{enumerate}
\item $\|\psi_n\|_2=1$.
\item $\rho(\xi):= 
\sum\limits_{j\in \Z} |\hat \psi_n (2^j\xi )|^2 = 1
\quad{\rm for~a.e.}~\xi\in\R.$
\item $t_q(\xi):=
\sum\limits_{j\geq 0} \hat \psi_n (2^j\xi)\overline
    {\hat \psi_n \left(2^j(\xi +2q\pi)\right)}
= 0 \quad{\rm for~a.e.~}~\xi\in\R$~{\rm and~for~all}
~$q\in 2\Z + 1.$
\end{enumerate} 
{\it Proof of} (1):
\vskip 2mm
We have
\begin{eqnarray*}
\|\hat\psi_n\|_2^2 
& = & (b_n-a_n)+\Bigl(\frac{a_n}{2}+2^{n-1}\pi -c_n\Bigr)
+d_n-\Bigl(\frac{e_n}{2}+2^{n-1}\pi\Bigr) \\
&   &  + \frac{1}{2}\Bigl\{\left(\frac{e_n}{2}-\frac{a_n}{2}\right)
       +(e_n-a_n)+\Bigl(\frac{e_n}{2}-\frac{a_n}{2}\Bigr) 
       +(e_n-a_n)\Bigr\} \\
& = & b_n-a_n+d_n-c_n+e_n-a_n = 2\pi.
\end{eqnarray*}
Therefore, $\|\psi_n\|_2 =1$. \pagebreak

\noindent{\it Proof of} (2):
\vskip 2mm

Let $\xi>0$. Since $\rho(\xi )= \rho(2\xi )$, it is enough 
to show that $\rho(\xi) = 1$ for a.e. 
$\xi\in [\alpha,2\alpha]$ for some $\alpha >0$. We will 
prove that $\rho(\xi)=1$  for a.e. 
$\xi\in [a_n,2a_n]=[a_n,b_n]=[a_n,e_n]\cup [e_n,b_n]$.

Suppose $\xi\in [a_n,e_n]$. Then $2^j\xi\in F_n$ 
only when $j=0,-1$. So,
$\rho(\xi )
=|\hat\psi (\xi )|^2+|\hat\psi (\frac{\xi}{2})|^2 
= \left(\frac{1}{\sqrt2}\right)^2 
            +\left(\frac{1}{\sqrt 2}\right)^2 
= 1$.
Now, 
\begin{eqnarray*}
\xi\in [e_n,b_n]\Longleftrightarrow 2^{n-1}\xi 
& \in & [c_n,d_n]\\
& = & \left[c_n, \frac{a_n}{2}+2^{n-1}\pi\right]
\cup \left[\frac{a_n}{2}+2^{n-1}\pi, 
      \frac{e_n}{2}+2^{n-1}\pi\right] \\
& & \cup \left[\frac{e_n}{2}+2^{n-1}\pi,d_n\right] \\
& =& I_1\cup I_2\cup I_3, ~{\rm say}.
\end{eqnarray*}
If $2^{n-1}\xi\in (I_1\cup I_3)$, then 
$2^j(2^{n-1}\xi)\in F_n$ only when $j=0$ 
(see the table). So, 
$\rho(\xi )=|\hat\psi_n (2^{n-1}\xi )|^2=1$. 
Also if $2^{n-1}\xi\in I_2$, then $2^j(2^{n-1}\xi)\in F_n$ 
only when $j=0 $ or $1$. 
So $\rho(\xi )
=|\hat\psi_n (2^{n-1}\xi)|^2 +|\hat\psi_n (2^n\xi )|^2 
= (\frac{1}{\sqrt 2})^2 +(-\frac{1}{\sqrt 2})^2 =1$.
Hence, $\rho(\xi ) =1$ for a.e. $\xi>0$.

For $\xi<0,$ it suffices to show that $\rho(\xi )=1$ 
on $[-b_n,-a_n]$. In fact, from the table it is clear 
that if $\xi\in [-b_n,-a_n]$, then $2^j \xi\in F_n$ 
only when $j=0 $. So, $\rho(\xi )=|\hat\psi_n (\xi )|^2 =1$ 
for a.e. $\xi\in[-b_n,-a_n]$. Thus, we have proved that 
$\rho(\xi )=1$ for a.e. $\xi\in\R$. \\

\noindent{\it Proof of} (3):
\vskip 2mm

We now have to prove that
\[ 
t_q(\xi) = \sum\limits_{j\geq 0} \hat \psi_n (2^j\xi)
\overline{\hat \psi_n \left(2^j(\xi +2q\pi)\right)}
         = 0 \quad{\rm for~a.e.}~\xi\in\R~
{\rm and~for~all}~q\in 2\Z + 1. 
\]

Since $t_q(\xi) = \overline{t_{-q}(\xi+2q\pi)}$ 
(see (\ref{eqn:tm}), it is enough to show that $t_q=0$ a.e., 
if $q$ is a negative odd integer. Suppose $q\not=-1,$ and is odd. 
We have $2^jq\neq 0,\pm 2^{n-1}\pm 2^{n-2}$. Therefore, 
if $2^j\xi\in F_n$, then from the table we observe that 
$2^j\xi+2\cdot 2^jq\pi\not\in F_n$, which shows that each 
term of the sum $t_q(\xi)$ is $0$. Hence, $t_q=0$ a.e.

It remains to prove that $t_{-1}(\xi)=0$ for a.e. $\xi\in\R$.
\[ 
t_{-1}(\xi)=\sum\limits_{j\geq 0} \hat \psi_n (2^j\xi)
\overline{\hat \psi_n (2^j\xi -2\cdot 2^j\pi)}. 
\]
From the table (see the 8th and 6th row), we observe 
that both $2^j\xi$ and $2^j\xi-2\cdot 2^j\pi$ belong to 
$F_n$ only in the following two cases :
\begin{itemize}
\item[(i)]  $2^j\xi \in [a_n+2^n\pi,e_n+2^n\pi]$ and $ j=n-1.$ 
\item[(ii)] $2^j\xi \in [\frac{a_n}{2}+2^{n-1}\pi,
             \frac{e_n}{2}+2^{n-1}\pi]$ and $ j=n-2.$ 
\end{itemize}
But both are equivalent to saying that 
$2^{n-1}\xi \in [a_n+2^n\pi,e_n+2^n\pi]$.
So we get $ t_{-1}(\xi)=0$, 
if $2^{n-1}\xi \not\in [a_n+2^n\pi,e_n+2^n\pi]$.

Now if  $2^{n-1}\xi \in [a_n+2^n\pi,e_n+2^n\pi]$, then
\begin{eqnarray*}
t_{-1}(\xi) & = & \hat\psi_n(2^{n-2}\xi)
\overline{\hat\psi_n(2^{n-2}\xi-2^{n-1}\pi)}
                 +\hat\psi_n(2^{n-1}\xi)
\overline{\hat\psi_n(2^{n-1}\xi-2^n\pi) }\\
            & = & \tfrac{1}{\sqrt2}\cdot \tfrac{1}{\sqrt2}
+\bigl(-\tfrac{1}{\sqrt2}\bigr)
                  \cdot\tfrac{1}{\sqrt2} \\
            & = & 0.
\end{eqnarray*}
We have proved that $t_q(\xi)=0$ for a.e. $\xi\in\R$ 
and for all $ q\in 2\Z + 1$.
Therefore, $\psi_n$ is a wavelet.

Our claim now is that $\psi_n\in\M_{n-2}$. For this, 
we have to show that
\begin{itemize}
\item[(a)] $F_n\cap(F_n+2^k q\pi) = \emptyset$ for all 
           $q\in 2\Z+1$ and $k=1,2,\dots n-2 $; and
\item[(b)] there exists a subset $E$ of $F_n$ such that
           $E+2^{n-1} q\pi\subset F_n$ for some $q\in 2\Z+1$.	
\end{itemize}
As $2^k q\pi = 2\cdot 2^{k-1} q\pi$, by referring to the 
table we observe that if $\xi\in F_n$, then 
$\xi+2\cdot 2^{k-1} q\pi\not\in F_n $ (since $q$ is odd 
and $k-1\leq n-3)$. So (a) is proved. To prove (b) 
consider $E = [\frac{a_n}{2},\frac{e_n}{2}]$ and $q=1$. 
Then, 
$E+2^{n-1}q\pi=[\frac{a_n}{2},\frac{e_n}{2}]+ 2^{n-1}\pi\subset F_n$. 

Since $n\geq 3$, we have proved that the equivalence classes 
$M_n,n\geq 1$ are non-empty.
\qed

\subsection {A family of wavelets belonging to the class $\M_0$}
\label{subsec:equiv2}
For $n\geq 3$, let $a_n,b_n,c_n,d_n$ and $e_n$ be as above.
Define the function $b$ on $[e_n,b_n]$ as follows. 
\begin{eqnarray*}
b(\xi)= 
\begin{cases}
\frac{1}{\sqrt 2} & {\rm if}~~\xi\in[e_n,\pi]\\
0                 & {\rm if}~~\xi\in [\pi,b_n].\\
\end{cases}
\end{eqnarray*}
Then we extend $b$ to the whole of $S_n$ by using (i)--(iv) 
of Theorem \ref{thm:sn}. Observe that 
$[e_n,b_n]\cap {\rm supp}~b\cap\bigl(\frac{1}{2^{n-1}}
{\rm supp}~b\bigr)=[e_n,\pi]$. We define $\theta$ on $\R$ as
\begin{eqnarray*}
\theta(\xi)=
\begin{cases}
\pi & {\rm if}~~\xi\in[e_n,\pi]\\
0                & {\rm otherwise.}\\
\end{cases}
\end{eqnarray*}
Clearly, $\theta$ satisfies the functional equation in (v) 
of Theorem \ref{thm:sn}. This choice of $b$ and $\theta$ 
will give us the wavelet $w_n$, where
$\widehat{w}_n(\xi)=e^{i\theta(\xi)}b(\xi)$. That is,
\begin{eqnarray*}
\widehat{w}_n(\xi)=
\begin{cases}
1                & {\rm if}~~\xi\in[-\pi,-a_n]\cup[a_n,e_n]
\cup[2^{n-1}\pi,d_n]\\
\frac{1}{\sqrt 2} & {\rm if}~~\xi\in[-d_n,-2^{n-1}\pi]
\cup[-b_n,-\pi]\cup[c_n,2^{n-1}\pi]\\
-\frac{1}{\sqrt 2} & {\rm if}~~\xi\in[e_n,\pi]\\
0                & {\rm otherwise.}\\
\end{cases}
\end{eqnarray*}
Weber \cite[Theorem 4]{web} also proved that a wavelet 
$\psi\in\L_1$ only if supp $\hat\psi$ is not partially 
self-similar with respect to any odd multiple of $2\pi$. 
In view of this result, to show that 
$w_n\in\M_0=\L_0\setminus\L_1$, we only have to find 
a set $H_n$ in supp $\widehat{w}_n$ such that $H_n+2q\pi$ 
is also a subset of supp $\widehat{w}_n$, for some odd 
integer $q$. The choice $H_n=[-b_n,-\pi]$ and 
$q=1$ will do the job, as 
$[-b_n,-\pi]+2\pi= [e_n,\pi]\subset {\rm supp}~\widehat{w}_n$.
\medskip
\begin{remark}
We can easily see that the wavelets constructed in this 
section are non-MRA wavelets.

(i) For wavelets of \S\ref{subsec:equiv1}: 
Let $\xi\in[\tfrac{a_n}{4},\tfrac{e_n}{4}]$. Then, 
$2\xi\in[\tfrac{a_n}{2},\tfrac{e_n}{2}]$,
$2(\xi+2\cdot 2^{n-3}\pi)\in[\tfrac{a_n}{2}+2^{n-1}\pi,
\tfrac{e_n}{2}+2^{n-1}\pi]$, 
$2^2\xi\in[a_n,e_n]$, and $2^2(\xi+2\cdot 2^{n-3}\pi)
\in[a_n+2^n\pi,e_n+2^n\pi]$. Hence, 
\bes
D_{\psi_n}(\xi) & \geq & |\hat\psi(2\xi)|^2 
+|\hat\psi(2(\xi+2\cdot 2^{n-3}\pi))|^2 
+|\hat\psi(2^2\xi)|^2 +|\hat\psi(2^2(\xi+2\cdot 2^{n-3}\pi))|^2 \\
                &  =   & \tfrac{1}{2}+\tfrac{1}{2}
+\tfrac{1}{2}+\tfrac{1}{2} = 2.                     
\ees
Therefore, $\psi_n$ is a non-MRA wavelet. 

(ii) For wavelets of  \S\ref{subsec:equiv2}:
Since supp $\widehat{w}_n\subseteq S_n$, the wavelets 
$w_n$ are non-MRA wavelets, 
by Theorem \ref{thm:sn}.
\end{remark}

\bibliographystyle{amsplain}

\end{document}